\documentclass[10pt]{amsproc}
\usepackage{amsmath,amssymb, graphicx, amscd,latexsym,color,comment}

\makeatletter
\newtheorem{Theorem}{Theorem}[section]
\newtheorem{Lemma}[Theorem]{Lemma}

\theoremstyle{definition}

\theoremstyle{remark}
\newtheorem{Remark}[Theorem]{Remark}

\newcommand{\eps}{\varepsilon}

\newcommand\vphi{\varphi}

\newcommand\si{\sigma}

\newcommand\de{\delta}

\newcommand\BC{ {\mathbb C}}

\newcommand\bfw{\mbox {\bf  w}}

\newcommand\bfz{\mbox {\bf  z}}

\newcommand\nl{\newline}

\newcommand\id{\rm{id}}

\newcommand\codim{\rm{codim}\/}

\newcommand\inv{^{-1}}

\def\inv{^{-1}}

\begin{document}

\title[On the connectivity of Milnor fiber for mixed functions
]
{On the connectivity of Milnor fiber for mixed functions
 }

\author
[M. Oka ]
{Mutsuo Oka\\\\
Dedicated to Professor L\^e  D$\tilde u$ng Tr\'ang for his 70's birthday}
\address{
Department of Mathematics,
Tokyo  University of Science、
Kagurazaka 1-3, \newline
Shinjuku-ku,
Tokyo 162-8601}
\email{oka@rs.kagu.tus.ac.jp}
\date{April 21,2018}
\begin{abstract}
In this note, we prove the connectivity of the Milnor fiber  for a mixed polynomial $f(\mathbf z,\bar{\mathbf z})$, assuming the existence of a sequence of smooth points of $f\inv(0)$ converging to the origin. This  result gives also a another proof for the connectivity of
the Milnor fiber of a non-reduced complex analytic function which is proved by A. Dimca
\end{abstract}
\maketitle

\maketitle

  \section{Introduction}
Let $f(\bfz,\bar\bfz)=\sum_{\nu,\mu}c_{\nu,\mu}\bfz^{\nu}{\bar\bfz}^\mu$ be a mixed polynomial of $n$-variables  $\bfz=$
\nl $(z_1,\dots, z_n)
\in \BC^n$ which has a singularity at the origin.
We say that $f$ satisfies {\em    Hamm-L\^e condition} at the origin if 
there exists a positive number $r_0$ such that for any $0<r_1\le r_0$, there exists a positive number $\delta(r_1)$ such that 
 the hypersurface $f\inv(\eta)\cap B_{r_0}^{2n}$ is non-empty,  non-singular and it intersects transversely with the sphere $S_r^{2n-1}$ for any $r$, $\eta$ with
$r_1\le r\le r_0$ and $0\ne |\eta|\le \delta(r_1)$. 
Here $B_r^{2n}$ is the ball $\{\mathbf z\in \mathbb C^n\,|\,\|\mathbf z\|\le r\}$ of radius $r$ and $S_r^{2n-1}$ is the boundary sphere of $B_r^{2n}$.
We call such a positive number  $r_0$ a stable Milnor radius.
If $f$ is a holomorphic function, 
 Hamm-L\^e condition  is always satisfied (Hamm-L\^e\,\cite{Hamm-Le}).
 For a strongly non-degenerate mixed function, this condition is also satisfied provided $f(\mathbf z,\bar{\mathbf z})$ is  either convenient (\cite{OkaMix}) or  locally uniformly  tame along vanishing coordinate subspaces (\cite{OkaAf}).
The following assertion is immediate from  Hamm-L\^e condition and Ehresman's fibration theorem (\cite{W}).
\begin{Lemma}[\cite{Hamm-Le,OkaAf}]
Assume that $f(\bfz,\bar{\bfz})$ be a mixed function of $n$-variables with $n\ge 2$ which satisfies the   Hamm-L\^e condition.
Choose $r_0>0$ as above and take  arbitrary positive numbers $r,r_1$, with $0<r_1\le r\le r_0$ and take a positive number $\de$ with $\de\le \de(r_1)$.
Consider the following tubular set and its boundary
\[\begin{split}
&E(r,\delta):=\{\mathbf z\in B_r^{2n} |\,|f(\mathbf z)|\le \delta\}\\
&\partial E(r,\delta):=\{\mathbf z\in B_r^{2n}\,|\, |f(\mathbf z,\bar {\mathbf z})|=\delta\}.
\end{split}
\]
Then  the mappings
\[\begin{split}
&f:E(r,\delta)\setminus f^{-1}(0)\to D_\delta\setminus\{0\}\,\,\text{and}\\
&f:\partial E(r,\delta)\to S^1_{\delta}
\end{split}
\]
are  locally trivial fibrations where $D_\de:=\{\zeta\,|\, |\zeta|\le \de\}$ and $S_\de^1=\partial D_\de$ (the boundary circle). They are homotopically equivalent  and their isomorphism classes do not depend on the choice of $r$ and $\delta$.
\end{Lemma}
This is called  {\em the tubular Milnor fibration} of $f$ at the origin and the fiber $f^{-1}(\de)\cap B_r$ is called the Milnor fiber at the origin.
It is known that the tubular Milnor fibration  is equivalent to the spherical Milnor fibration ( see\,\cite{Milnor,OkaMix})
\[
\vphi: S_{r}\setminus K\to S^1,\quad \vphi(z)=f(\bf z,\bar{\bf z})/|f(\bf z,\bar{\bf z})|
\]
\section{Statement of the result} 
\subsection{Connectivity of Milnor fibers of holomorphic functions} We first recall basic facts  about the connectivity of the Milnor fibers of  holomorphic functions.
Assume that $f: (U,{\bf 0})\to (\mathbf C,0) $ be a holomorphic function with a singularity at the origin
where $U$ is an open neighborhood of the origin.
A fundamental result for the connectivity is
\begin{Theorem}[Milnor \cite{Milnor}]
 The Milnor fiber $F$ has a  homotopy type of an $(n-1)$-dimensional CW-complex. 
If further the origin is an isolated singularity,
$F$ is $(n-2)$-connected.
\end{Theorem}
Let $\vphi:S_r^{2n-1}\setminus K\to S^1$ be the spherical Milnor fibration.
The proof depends on Morse theory. 
Milnor proved that the fiber has a homotopy type of an  $n-1$ dimensional CW-complex, by showing that the index of a suitable Morse function is less than or equal to $n-1$. Secondly, if the origin is an isolated singularity at the origin, $\bar F$ is a manifold with boundary and the inclusion $F\subset \bar F$ is a homotopy equivalence 
in  the sphere $S_r^{2n-1}$ and   the vanishing $H_j(F)=0$ for $0<j<n-1$ follows from
the Alexander duality and the homotopy equivalence $F_{\pi}\subset S_r^{2n-1}\setminus \bar F$
where $F_{\pi}=\vphi\inv(-1)$. The simply connectedness of $F$ for $n\ge 3$ follows from handle body argument.

 Kato and Matsumoto further generalized this assertion as follows.
\begin{Theorem}[Kato and Matsumoto  \cite {Kato-Matsumoto}] \label{Kato-Matsumoto} Assume that $s$ is the  dimension of the critical points locus at the origin.
Then  $F$ is  $(n-s-2)$-connected
\end{Theorem}
Their proof depends on  the above  result of Milnor and an inductive argument using Whitney stratification.
\subsection{Mixed functions}
 For mixed functions, 
 no similar connectivity statement is  known. 
  A main reason is that  the tangent space of a mixed hypersurface has no complex structure
   and  Morse function argument does not help so much as in the holomorphic case.

The following is our result which is a first  step for the connectivity
of the Milnor fiber of a mixed function. For the proof, we use  an elementary but completely different viewpoint which is nothing to do with  Morse functions. We use a one-parameter group of diffeomorphisms
which containes the monodromy map. 
\begin{Theorem}\label{main-theorem}
Assume that $f(\mathbf z,\bar{\mathbf z})$ is a mixed function of $n$-complex variables
$\mathbf z=(z_1,\dots,z_n)$ with $n\ge 2$
which satisfies Hamm-L\^e condition at the origin and take  a stable Milnor radius $r_0$
and
assume that there exists a mixed smooth point  $\mathbf w\in f\inv(0)$ with $\|\mathbf w\| < r_0$
and the sphere of radius $\|\mathbf w\|$ intersect $ f\inv(0)$ transversely at $\mathbf w$. 
We assume also that $\dim_{\mathbb R}f^{-1}(0)=2n-2$.
Then the Milnor fiber is connected.
\end{Theorem}
\begin{Remark}Suppose that $f$ is a holomorphic function and   $s$ is the dimension of the critical points locus of $f$ at the origin.
 If $f$ is non-reduced, then  $s=n-1$ and  Theorem \ref{Kato-Matsumoto} says nothing about the connectivity.
\end{Remark}
\subsection{Proof of Theorem \ref{main-theorem}}Let $r_0$ be a positive number which satisfies the
Hamm-L\^e condition.
Take a mixed  smooth point $\bfw\in f^{-1}(0)$ with $\|\bfw\|< r_0$. Put $r:=\|\bfw\|$.
Fix  positive numbers $r_1$ with $r_1<r$ and  $\delta\le \delta(r_1)$
and we consider the Milnor fibration
\[
(\star)\quad
f:\partial E(r_0,\delta)\to S_\delta^1.
\]
Let $F_{\theta}:=f^{-1}(\delta e^{i\theta})\cap B_{r_0}$ for $\theta\in \mathbb R$ be the Milnor fiber 
over $\delta e^{i\theta}$.
We assume also the sphere with radius $r$, $S_r^{2n-1}$ intersects transversely with the hypersurface
$f^{-1}(0)$ at the smooth point $\bfw$. For the proof, we use a certain one-parameter family of diffeomorphisms associated with the tubular Milnor fibration ($\star$).
Recall that a one-parameter family of the characteristic diffeomorphisms $h_\theta,\,\theta\in \mathbb R$ are constructed by integrating a given horizontal vector field $\mathcal V$. Here a vector field $\mathcal V$ is  called {\em a  horizontal
vector field } on 
$\partial E(r_0,\delta)$ if it satisfies  the following property:
\[
Tf_{p}(\mathcal V(p))=\frac{\partial}{\partial \theta}(f(p)),\quad p\in \partial E(r_0,\de)
\]
where  $Tf_p:\,T_p\partial E(r_0,\delta)\to T_{f(p)}S_\delta^1$ is the tangential map and furthermore $\mathcal V$ is
tangent to $S_{r_0}\cap\partial E(r_0,\delta)$ on the boundary and 
 $\frac{\partial}{\partial \theta}$ is the unit angular vector field along  $S_{\de}^1$. In other word, 
$\frac{\partial}{\partial \theta}$ is the tangent vector of the curve $t\mapsto \de e^{it}$.
As $\partial E(r_0,\de)$ is a compact manifold with boundary, there is an integral
$\vphi:\partial E(r_0,\de)\times \mathbb R\to \partial E(r_0,\de)$
such that 
$\vphi(p,0)=p$ and $\vphi(p,t),\,-\infty<t<\infty$ is the integral of $\mathcal V$ starting at $p$ for $t=0$.
Let $h_{\theta}:\partial E(r_0,\delta)\to \partial E(r_0,\delta),\,\theta\in \mathbb R$ be the corresponding one parameter family of  characteristic diffeomorphisms, which are defined by
$h_\theta(p)=\vphi(p,\theta)$. Note that 
$\{h_\theta\}$ satisfy the property
\[
h_\theta(F_\eta)= F_{\eta+\theta},\quad \eta,\theta\in \mathbf R.
\]
They also satisfy the equalities
\[
h_0={\id},\quad
h_\theta\circ h_{\xi}=h_{\theta+\xi},\quad \theta,\xi\in \mathbb R. 
\]
Using one-parameter family $\{h_\theta,\theta\in \mathbb R\}$, the monodromy map $h:F\to F$ is given by the restriction $h_{2\pi}|F$ with $F=F_0$.

Now we construct $\mathcal V$ more carefully.
We  take a local real-analytic  chart $U$ with local coordinates $(x_1,y_1,\dots, x_n,y_n)$
centered at the given smooth point  $\bfw\in f\inv(0)$ such that $ f(x_1,y_1,\dots, x_n,y_n)=x_n+iy_n$. We take $\de$ sufficiently small and 
we take a normal disc $D$ of $f^{-1}(0)$ centered at $\bfw$
so that $|f(\bfz)|=\delta$ for any  $\bfz\in \partial D$ and thus $\partial D\subset \partial E(r_0,\delta)$.
In the above coordinates, we can assume $D=\{(0,\dots, 0,x_n,y_n)\,|\,x_n^2+y_n^2\le \de^2\}$.
By the transversality assumption, we may assume that $\mathcal V$ is tangent to $\partial D$ which implies that
$h_{2\pi}(\bfz)=\bfz$ for $\mathbf z\in \partial D$.
%
Note that the family of tubular neighborhoods
$\{E(r,\delta)|\delta< \de(r)\}$ and the family of disk neighborhoods $\{B_s^{2n},s>0\}$ are cofinal neighborhood systems of the origin. 
Thus 
 $E(r,\delta)$ is  contractible and $E\setminus f^{-1}(0)$ is connected by the assumption  $\codim_{\mathbb R}\,f^{-1}(0)=2$, $\partial E(r,\delta)$ is also connected,
as $E(r_0,\delta)\setminus f^{-1}(0)$ and $\partial E(r_0,\delta)$ are homotopy equivalent.

Now we are ready to prove the connectivity of the Milnor fiber $F_0$.
Fix a point $p\in \partial D\cap F_0$ and take an arbitrary point $q\in F_0$.
As  $\partial E(r,\delta)$ is  connected,
 we can find a path
$\si:[0,1]\to \partial E(r_0,\delta)$ such that $\si(0)=q$ and $\si(1)=p$.
 We will show that $q\in F_0$ can be joined
 to $p$ by a path in the fiber $F_0$, modifying the path $\si$. 
 Consider first the closed loop $f\circ \si:(I,\{0,1\})\to 
(S^1_{\delta},\delta)$ and let $m$ be the rotation number, that is
\[
m=\frac 1{2\pi}\int_{0}^1 d\arg\,f(\sigma(\theta)).
\]
Let $\omega:I\to \partial D$ be the clockwise rotation $m$-times along the boundary of $D$
starting at $p$. 
 Consider the path $\sigma'$ which is given as  the composition of paths $\si':=\si\cdot\omega:I\to \partial E(r,\delta)$
 and define a function $\psi:\mathbb [0,1]\to \mathbb R$ by
 \[
\psi(t)=\int_0^t d\arg f(\si'(t)) dt.
\]
Note that
 $\psi(t)\equiv \arg f(\si'(t))$ modulo $2\pi$ and $\psi(0)=\psi(1)=0$. This follows from
the observation that the rotation number of $f\circ \si'$ is zero by the definition of $\omega$. 
Now we can deform $\si'$  using  characteristic diffeomorphisms $h_t,\,t\in \mathbb R$  into a path in the fiber $F_0$ as follows. Define a modified path
\[
\hat{\si}(t):=h_{-\psi(t)}(\si'(t)),\quad 0\le t\le 1.
\]
Then we have the equality
\[
\arg\,f(\hat{\si}(t))=-\psi(t)+\arg f(\si'(t))\equiv 0,\quad 0\le t\le 1.
\]
 Thus the deformed path $\hat \sigma$ is entirely 
included in $F_0$. 
By the construction,  $\hat\si(0)=q$ and $\hat\si(1)=p$  and $\hat \sigma$ is the path which connect
$q$ and $p$ in $F_0$. \qed
\subsection{Generalization of Main Theorem for holomorphic functions}
The following is known by A. Dimca \cite{Dimca}.
We give another direct proof of this assertion, using a similar argument as in the proof of Theorem \ref{main-theorem}.
\begin{Theorem}{\rm ( Prpoposition 2.3, \cite{Dimca} )} Assume that $f$ is a germ of holomorphic function at the origin of $\mathbb C^n,\,n\ge 2$
with $f(\mathbf 0)=0$ and assume that
$f$ is factored as $f_1^{n_1}\cdots f_r^{n_r}$ where $f_1,\dots,f_r$ are  irreducible in $\mathcal O_n$
and mutually coprime.
Put $n_0=\gcd(n_1,\dots, n_r)$.  Then Milnor fiber of $f$ at the origin is 
connected if and only if $n_0=1$.
\end{Theorem}

Here $\mathcal O_n$ is  the ring of germs of holomorphic functions at the origin. 
\begin{proof}Let $F$ be the Milnor fiber of $f$.
We can write $f=g^{n_0}$  with
 $g=f_1^{n_1/n_0}\cdots f_r^{n_r/n_0}$
and $F$ is diffeomorphic to  disjoint sum of  $n_0$ copies of  the Milnor fibers of $g$.
Thus $F$ is not connected if $n_0\ge 2$.
Assume $n_0=1$.
Take a stable radius $r_0$ which satisfies Hamm-L\^e condition.
For each $1\le i\le r$, we consider the reduced irreducible component $V_i=\{f_i=0\}$ and take a non-singular point $p_i\in V_i\setminus \bigcup_{j\ne i}V_j$ with $\|p_i\| < r_0$.
 Choose a positive number  $r<r_0$ with $r\le \|p_i\|, i=1,\dots,r$ and we consider  Milnor fibration 
\[
f:\partial E(r_0,\de)\to S_\de^1,\quad \de\le \de(r) 
\]
We assume  that in a sufficiently small neighborhood  $U_i$ of $p_i$, any analytic branch of $f^{1/n_i}$ is  a well-defined
 single-valued function.  Taking a one branch $\widetilde{f_i}$ and we take 
$\widetilde{f_i}$  as the last coordinate function of  a complex  analytic coordinate system $\mathbf z_i=(z_{i1},\dots, z_{in})$ 
in $U_i$ i.e., $\widetilde{f_i}(\mathbf z_i)=z_{in}$. Consider a normal disk at $p_i$ which is defined by 
$D_i=\{(0,\dots, 0,z_{in})\,|\,|z_{in}|\le {\de}^{1/n_i}\}$
where we assume that $\de$ is sufficiently small so that $D_i\subset U_i$.
Note that $\widetilde{f_i}$ is locally written as $\widetilde{f_i}=f_i\cdot u_i$ where $u_i$ is a unit in $U_i$ and $f(\partial D_i)=S_\de^1$
and $\partial D_i\subset \partial E(r,\de)$.
Construct a one parameter family of chracteristic diffeomorphisms $h_t,\,t\in \mathbb R$ as before such that
$\partial D_i$ is stable by $h_t$ and $h_t|\partial D_i$ is the rotation of angle $t/n_i$, under the identification
$\widetilde{f_i}:\partial D_i\to S_{\de^{1/n_i}}^1$. (Note that $f:\partial D_i\to S_\de^1$ is the $n_i$ rotation.)
$F\cap\partial D_i$ can be identified with $n_i$-th roots of $\de$ and we put them as
$F\cap\partial D_i=\{p_{i,0},\dots,p_{i,n_i-1}\}$ so  that  the monodromy
$h:=h_{2\pi}$ acts simply as a cyclic permutation $p_{i,j}\mapsto p_{i,j+1}$. 
Using a similar discussion as in the proof of Theorem \ref{main-theorem}, for a given $q\in F$,
 we can connect $q$  to some point in $F\cap D_j$ for any $j$. To see this, we first take a path $\sigma$ so that $\si(0)=q$ and $\si(1)=p_{j,0}$, then 
in the argument of the proof of  Theorem \ref{main-theorem}, we simply replace $\omega$ by $m/n_j$ rotation in the clockwise direction along $\partial D_j$
and  do the same argument, where $m$ is the rotation number of $f\circ \sigma$. Note that $\omega$ need not a closed loop
but the image of $\omega$ by $f$ is a closed loop and it gives $-m$ rotation.
Thus the proof is reduced to show that the points  $\{p_{1,0},\dots, p_{1,n_1-1}\}$ are in the same connected component of $F$.

In particular taking $q=p_{1,0}$, we can find 
a path $\ell_j$ in $F$ which connects $p_{1,0}$ to some $p_{j,\nu_j}$  for any  $j,\,j=2,\dots,r$. 
Let $\ell_j^{a}=h^a(\ell_j)$ the image of $\ell_j$ by $h^{a}$. Then it connects $p_{1,a}$ to $p_{j,\nu_j+a}$ in $F$.
Now we consider the image of $\ell_j^{a},\,a=0,\dots,n_1-1$ under
$h^{n_j}$ which fixes $p_{j,\nu_j+a}$ but the other end $p_{1,a}$ of $\ell_j^{a}$ goes to $p_{1,n_j+a}$. 
As $F$ is stable by the monodromy map $h=h_{2\pi}$, this image is also a path in $F$.
Thus $p_{1,a}$ and $p_{1,n_j+a}=h^{n_j}(p_{1,a})$ are connected by the path
$\ell_j^{a}\cdot( h^{n_j}\circ (\ell_j^{a})^{-1})$ in $F$ for any $a$. See Figure 1 for $a=0$.
As we have assumed  $\gcd(n_1,\dots,n_r)=1$, there exist integers $a_1,\dots,a_r$ so that we can write $1=\sum_{i}^r a_in_i$. 
Then 
\[
p_{1,1}=h(p_{i,0})=h^{\sum_{i=1}^r a_i n_i}(p_{1,0})=h^{a_r n_r}(\cdots(h^{a_2n_2}(p_{1,0})\cdots)).
\]
Put $p_{1,\mu_2}:=h^{a_2n_2}(p_{1,0})$ and $p_{1,\mu_{j+1}}:=h^{n_{j+1}a_{j+1}}(p_{1,\mu_j})$
for $j=2,\dots, r-1$.
Then points $p_{1,\mu_2},\dots,p_{1,\mu_r}=p_{1,1}$ are in the same connected component of $p_{1,0}$
in $F$. Thus $p_{1,1}$ and $p_{1,0}$ are in the same component.
Repeating the same argument, we conclude that $p_{1,0},\dots, p_{1,n_1-1}$ are all in the same connected component of $F$. Combining the above observation, this proves that $F$ is connected.
\begin{figure}[htb]
\setlength{\unitlength}{1bp}
\begin{picture}(600,250)(-70,-640) 
{\includegraphics{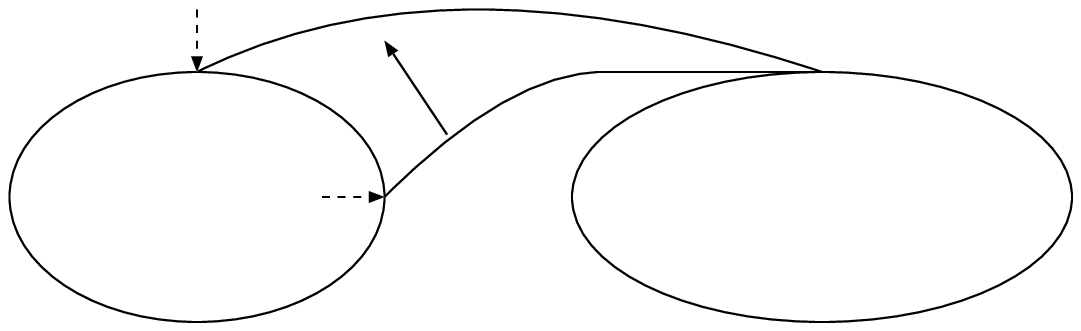}}
\put(-290,-540){$p_{1,0}$}
\put(-320,-480){$p_{1,n_j}$}
\put(-260,-510){$h^{n_j}$}
\put(-230,-530){$\ell_j^0$}
\put(-250,-470){$h^{n_j}(\ell_j^0)$}
\put(-310,-530){$D_1$}
\put(-130,-530){$D_j$}
\put(-130,-490){$p_{j,\nu_j}$}
\end{picture}
\vspace{-2cm}
\caption{$\ell_{1,0}$ and $h^{n_j}(\ell_{1,0})$}\label{ell-j}
\end{figure}

\end{proof}
\begin{Remark}
The special case $f=z_1^{a_1}\cdots z_k^{a_k}$ has been considered in  Example (3.7), \cite{Okabook}.
For mixed function case, 
the above argument does not work.
There does not exist any correspondence between irreducible components  of $f\inv(0)$ as real algebraic sets and the irreducible factors of $f(\mathbf z,\bar{\mathbf z})$ as a polynomial in $\mathbb C [\mathbf z,\bar{\mathbf z}]$
or a germ of an analytic function of $2n$ variables. For example, let $f(z,w,\bar z,\bar w)$ be the strongly homogeneous mixed polynomial of two variables which come from the Rhie's Lens  equation $\vphi_n=0$,
\[
\vphi_n(z):=\bar z-\frac{z^{n-2}}{z^{n-1}-a^{n-1}}-\frac{\eps}{z},\,n\ge 2
\]
where $\eps, a$ are sufficiently small positive numbers $0<\eps\ll a\ll 1$.
$f$ is obtained as the homogenization of the numerator of $\vphi_n$. Then $f$ is irreducible but $f\inv(0)$ has $5n-5$ line components in $\mathbb C^2$(\cite{OkaLens,Rhie}).
\end{Remark}
\section{Application}
\subsection{Convenient strongly non-degenerate mixed functions}
Let $f(\mathbf z,\bar{\mathbf z})$ be a strongly non-degenerate convenient mixed polynomial. Then 
 the Hamm-L\^e condition is  satisfied (Lemma 28, \cite{OkaMix}). Furthermore the mixed hypersurface  $V=f\inv(0)$ has an isolated 
 singularity at the origin (Corollary 20, \cite{OkaMix}). Thus the assumption of Theorem \ref{main-theorem}
 is satisfied and we have

\begin{Theorem}
Let $f(\bfz,\bar{\bfz})$ be a strongly non-degenerate convenient mixed function and $n\ge 2$.
Then the Milnor fiber is connected.
\end{Theorem}
\subsection{Non-convenient mixed polynomials}
Let $f(\mathbf z,\bar{\mathbf z})$ be
a strongly non-degenerate mixed polynomial which is not convenient.
In general, such a mixed hypersurface defined by $f$  has a non-isolated singular locus. Let $I$ be a non-empty subset of
$\{1,\dots, n\}$. 
We use the following notations.
\[\begin{split}
\mathbb C^I&=\{\mathbf z\,|\, z_j=0,\,\forall j\notin I\},\\
\mathbb C^{*I}&=\{\mathbf z\,|\, z_j=0 \iff j\notin I\}.
\end{split}
\]
We say $\mathbb C^I$ is {\em a vanishing coordinate subspace for $f$} 
if the restriction of $f$ to the coordinate subspace $\mathbb C^I$, denoted by $f^I$, is identically zero.
Otherwise, we say $\mathbb C^I$  {\em a non-vanishing coordinate subspace}.
Put $V^{\sharp}$ be the union of $V\cap {\mathbb C}^{*I}$ for all $I$ such that  $\mathbb C^I$ is a non-vanishing coordinate subspace. We know that $V^{\sharp}$ is mixed non-singular as a germ at the origin
(Theorem 19, \cite{OkaMix}).
 We say that {\em $f$  is uniformly locally tame on the vanishing coordinate subspace $\mathbb C^{I}$} if there exists a positive number $\eps>0$ and 
for any non-negative weight vector $P=(p_1,\dots, p_n)\in {\mathbb Z}_{\ge 0}^n$ with
$I(P)=I$ where $I(P):=\{i\,|\, p_i=0\}$, the face function $f_P$ is strongly non-degenerate as a function of 
$\{z_j\,|\, j\notin I\}$ for any fixed  $\mathbf z_I=(z_i)_{i\in I}\in \mathbb C^{*I}$ with $\sum_{i\in I}|z_i|^2\le \eps$.
A strongly non-degenerate mixed function which 
is uniformly locally tame along any vanishing coordinate subspaces satisfies 
 Hamm-L\^e condition ( Proposition 11, \cite{OkaAf}).
Thus  we have
\begin{Theorem}
Assume that $f(\bfz,\bar{\bfz})$ is  a strongly non-degenerate  mixed function which is uniformly locally tame
along vanishing coordinate subspaces and $n\ge 2$. We assume also that $V^{\sharp}$ is a non-empty germ of mixed variety. Then the Milnor fiber is connected.
\end{Theorem}
A similar but weaker result is proved for a strongly mixed homogeneous polynomial (\cite{OkaConj}).
As a next working problem, we also propose a conjecture  about the fundamental groups of the Milnor fiber.


\bibliographystyle{amsalpha}

\end{document}